\documentclass[12pt]{amsart}
\usepackage{epsfig}
\usepackage{amscd}
\usepackage[mathscr]{eucal}
\usepackage{amssymb}
\usepackage{amsxtra}
\usepackage{amsmath}
\usepackage[all]{xy}

\newtheorem{thm}{Theorem}

\begin{document}

\title{Homology of Lagrangian Submanifolds in Cotangent Bundles}
\date{\today}
\author{Lev Buhovsky}
\address{Lev Buhovsky, School of Mathematical Sciences, Tel-Aviv
  University, Ramat-Aviv, Tel-Aviv 69978, Israel}
\email{levbuh@post.tau.ac.il}
\maketitle

\section{Introduction and main results} \label{S:intro}

The present note deals with topological restrictions on Lagrangian
submanifolds in cotangent bundles. This subject has attracted a
great deal of attention since the beginning of symplectic topology
(see e.g. \cite{Gr, LS, P, V}). In this note we concentrate on
cotangent bundles of spheres and of Lens spaces. Denote by
$T^{*}(S^{2k+1})$ the cotangent bundle to the $(2k+1)$ -
dimensional sphere, endowed with its canonical symplectic
structure. Henceforth all Lagrangian submanifolds are assumed to
be \textit{compact, connected and with no boundary}. Our first
result is:
\begin{thm} \label{T:sphere}
   Let $L\subset T^{*}(S^{2k+1})$ be a Lagrangian submanifold with
   $H_{1}(L,\mathbb{Z})=0$. Then the cohomology of $L$, with
   $\mathbb{Z}_{2}$ - coefficients, equals the cohomology of the
   zero-section, namely $H^{*}(L,\mathbb{Z}_{2})\cong
   H^{*}(S^{2k+1},\mathbb{Z}_{2})$.
\end{thm}

A few remarks are in order before we continue. Clearly,
Theorem~\ref{T:sphere} is not empty, since we have the
zero-section (and its images under symplectic diffeomorphisms) as
Lagrangian submanifolds of $T^{*}(S^{2k+1})$.  However, beyond
these we do not know any other examples of Lagrangians $L\subset
T^{*}(S^{2k+1})$ with $H_{1}(L,\mathbb{Z})=0$. In fact, a long
standing (folkloric) conjecture asserts that the only exact
Lagrangian submanifolds in cotangent bundles are isotopic to the
zero-section. Theorem~\ref{T:sphere} can be viewed as some
supporting evidence to this conjecture.

Theorem~\ref{T:sphere} admits the following generalization to Lens
spaces.  Denote by $Lens_m^{2k+1} = S^{2k+1} / \mathbb{Z}_m$ the
$(2k+1)$-dimensional Lens space.

\begin{thm} \label{T:lens}
   Let $L \subset T^{*}(Lens_{m}^{2k+1})$ be a Lagrangian submanifold
   with $ \pi_{1}(L)= \mathbb{Z}_{m}$.  Then
   $H^{*}(L,\mathbb{Z}_{2})=H^{*}(Lens_{m}^{2k+1},\mathbb{Z}_{2})$.
\end{thm}

The phenomenon arising from Theorems~\ref{T:sphere}
and~\ref{T:lens} can be described as some kind of "homological
rigidity'' (or uniqueness). Namely, low-dimensional topological
invariants of Lagrangian submanifolds determine their {\em entire
homology}. First examples of this phenomenon were discovered by
Seidel~\cite{Se-1} and then by Biran~\cite{Bi-ICM,
Bi-NonIntersect, Bi-Homological}.

 Very recently\footnote{After the first draft of this paper has been written.}  Siedel ~\cite{Se-2}
has generalized, using other methods, Theorem~\ref{T:sphere}. His
result deals with compact connected Lagrangian submanifold of $ L
\subset M = T^{*}(S^{n}) $ (where $ n $ can be also even) with  $
H^{1}(L , \mathbb{Z}) = 0 $ and vanishing second Whitney class $
\textit{w}_{2}(L) $. Under these assumptions, he proves (among
other results) uniqueness of cohomology with coefficients in the
field $ \mathbb{C} $, that is, $ H^{*}(L,\mathbb{C}) \cong
H^{*}(S^{n},\mathbb{C}) $. The proof of this results uses
different tools than ours.

 In comparison to Siedel's approach, our methods apply only to
odd-dimensional spheres, since we use in a crutial way existence
of a fixed point free circle action on the sphere . Also note that
our results are concerned with homology with coefficients in
$\mathbb{Z}_{2}$. However, we believe that essentially the same
proof should imply, with some corrections, homological uniqueness
with coefficients in $\mathbb{Z}$.

The rest of the paper is devoted to proving
Theorems~\ref{T:sphere} and~\ref{T:lens}. The proofs are based on
techniques of symplectic topology, the main ingredient being
computations in Floer homology. In Section~\ref{S:computations}
below we outline these computations and in Section~\ref{S:proofs}
we apply them in order to prove our theorems.

\section{Computation in Floer homology} \label{S:computations}
Before we prove our main theorems in Section~\ref{S:proofs} we
need to recall some important facts from Floer theory that will be
used in our proof. Most of the theory in this section is due to
Oh. We refer the reader to~\cite{Oh-Spectral} for more details
(see also~\cite{Se-1}, and~\cite{Bi-NonIntersect,
Bi-Homological}).

Let $(M, \omega)$ be a tame symplectic manifold (see~\cite{ALP}
for the definition), and let $L \subset (M, \omega)$ be a monotone
Lagrangian submanifold with minimal Maslov number $N_L \geq 2$
(see~\cite{Oh-HF1, Oh-Spectral}). For any Hamiltonially isotopic
copy of $L$, say $L'$, we denote by $HF(L,L')$ the Floer homology
of the pair $(L, L')$. Since this invariant is independent of the
choice of $L'$ we shall sometimes denote it also by $HF(L)$. We
now describe Oh's~\cite{Oh-Spectral} approach for computing
$HF(L)$ using a spectral sequence whose first step is the singular
cohomology $H^*(L;\mathbb{Z}_2)$ of $L$.  Throughout this paper we
work with Floer and Morse homologies with
$\mathbb{Z}_2$-coefficients.

Let $L_{\epsilon} \subset (M, \omega)$ be a small Hamiltonian
perturbation of $L$ defined using a Morse function $f:L \to
\mathbb{R}$ and a Weinstein tubular neighborhood $\mathcal{U}$ of $L$.
Then the total Floer complex $CF(L, L_{\epsilon})$ can be graded by
the Morse indices of $f$. We shall denote this grading by $CF^*(L,
L_{\epsilon})$. Then, the Floer differential $d_{F} : CF(L,
L_{\epsilon}) \to CF(L, L_{\epsilon})$ can be written as $d_{F} =
\partial_0 + \partial'$ where $\partial_0$ comes from counting Floer
trajectories lying entirely inside $\mathcal{U}$ while $\partial'$
comes from counting those trajectories going out of $\mathcal{U}$. By
Floer's work~\cite{F-1, F-2}, $\partial_0$ can be identified with the
differential of Morse homology of the function $f$, hence $\partial_0
: CF^*(L, L_{\epsilon}) \to CF^{*+1}(L, L_{\epsilon})$ (i.e.
$\partial_0$ raises grading by $1$) and $H^*(\partial_0) \cong H^*(L;
\mathbb{Z}_2)$. (Note however that the operator $\partial'$ is not a
differential.) It is shown in~\cite{Oh-Spectral} that the operator
$\partial'$ splits as $\partial' = \partial_1 + \dots +
\partial_{\nu}$, where $\nu = \big[ \frac{\dim L + 1}{N_L} \big]$, and
the behavior of each $\partial_i$ in terms of grading is $\partial_i :
CF^*(L, L_{\epsilon}) \to CF^{*+1 - i N_L}(L, L_{\epsilon})$. By a
straightforward algebraic computation (see~\cite{Oh-Spectral}) it
follows that $\partial_1$ descends to a homomorphism
$[\partial_1]:H^*(\partial_0) \to H^{*+1 - N_L}(\partial_0)$ with
$[\partial_1]^2 = 0$. Thus $[\partial_1]$ is a differential on the
total homology $\oplus_{k=0}^{\dim L} H^k(\partial_0)$. (Note that
$[\partial_1]$ is not compatible with the grading of $H^*(\partial_0)$
in the sense that it does not raise grading by $1$.)  Denote by
$H^*([\partial_1])$ the homology of $[\partial_1]$, namely
$$H^{k}([\partial_1]) = \frac{\ker \left([\partial_1]:H^i(\partial_0)
     \to
     H^{k+1-N_L}(\partial_0)\right)}{\textnormal{im\,}\big([\partial_1]:H^{k-1+N_L}(\partial_0)
  \to H^k(\partial_0)\big)} \quad \textnormal{for every } k.$$
Here
for the grading of $H^*([\partial_1])$ we use the one induced from
that of $H^*(\partial_0)$.

Continuing by induction we get that $\partial_i$ descends to a
homomorphism $[\partial_i]: H^*([\partial_{i-1}]) \to H^{*+1 - i
  N_L}([\partial_{i-1}])$ with $[\partial_i]^2 = 0$. We denote by
$H^*([\partial_i])$ the corresponding homology (with grading induced
from $H^*([\partial_{i-1}])$, namely
$$H^{k}([\partial_i]) = \frac{\ker
  \left([\partial_i]:H^k([\partial_{i-1}]) \to H^{k+1-i
       N_L}([\partial_{i-1}])\right)}
{\textnormal{im\,}\big([\partial_i]:H^{k-1 + i N_L}([\partial_{i-1}])
  \to H^k([\partial_{i-1}])\big)} \quad \textnormal{for every } k.$$
Finally, it is shown in~\cite{Oh-Spectral} that
$$\bigoplus_{k=0}^{\dim L} H^k([\partial_{\nu}]) \cong \ker d_F /
\textnormal{im\,} d_F = HF(L, L_{\epsilon}) \cong HF(L).$$

\section{Proofs of theorems} \label{S:proofs}

%
%

\begin{proof}[Proof of Theorem~\ref{T:sphere}]
   Let us first describe the main ideas of the proof. The first
   observation is that a small neighborhood of the zero-section of
   $T^{*}(S^{2k+1})$ can be symplectically embedded into
   $\mathbb{C}P^{k}\times \mathbb{C}^{k+1}$ (such an embedding was
   first described by Audin, Lalonde and Polterovich [A-L-P]).
   Therefore, if $L\subset T^{*}(S^{2k+1})$ is a Lagrangian
   submanifold we can Lagrangianly embed it also into
   $\mathbb{C}P^{k}\times \mathbb{C}^{k+1}$. The benefit of this
   embedding is that in $\mathbb{C}P^{k}\times \mathbb{C}^{k+1}$, due
   to the $\mathbb{C}^{k+1}$-factor, any compact subset can be
   disjoint from itself via a Hamiltonian isotopy (e.g. by a linear
   translation). In particular, viewing $L$ as a Lagrangian in
   $\mathbb{C}P^{k}\times \mathbb{C}^{k+1}$ we have $HF(L,L)=0$.
   Comparing this vanishing of Floer homology with the spectral
   sequence computation from Section~\ref{S:computations} , we shall derive our
   restrictions on the first step of this sequence which is
   $H^{*}(L,\mathbb{Z}_{2})$.

   Let us turn now to the details of the proof. We first use a
   Lagrangian embedding $S^{2k+1}\hookrightarrow \mathbb{C}P^{k}\times
   \mathbb{C}^{k+1}$ due to Audin, Lalonde and Polterovich [A-L-P]. To
   describe this embedding denote by $i:S^{2k+1}\subset
   \mathbb{C}^{k+1}$ the standard inclusion as the unit sphere and by
   $\pi :S^{2k+1}\rightarrow \mathbb{C}P^{k}$ the Hopf fibration.
   Next, we endow $\mathbb{C}P^{k}$ with its standard symplectic
   structure, normalized so that the symplectic area of a projective
   line is $\pi$. A simple computation shows that
   $$S^{2k+1}\hookrightarrow \mathbb{C}P^{k}\times \mathbb{C}^{k+1},
   \quad z \mapsto (\pi(z),\overline{i(z)}).$$
   is a Lagrangian
   embedding (here, $\overline{(\cdot )}$ stands for usual complex
   conjugation).

   By Darboux-Weinstein theorem we now have a symplectic embedding of
   a small tubular neighborhood of the zero-section of
   $T^{*}(S^{2k+1})$ into $\mathbb{C}P^{k}\times \mathbb{C}^{k+1}$.
   Using homotheties along the cotangent fibers of $T^{*}(S^{2k+1})$
   we may assume that the Lagrangian $L$ lies in the above small
   tubular neighborhood of the zero-section. Thus we obtain a
   Lagrangian embedding of $L$ into $\mathbb{C}P^{k}\times
   \mathbb{C}^{k+1}$. From now on we shall view $L$ as a Lagrangian
   submanifold of $\mathbb{C}P^{k}\times \mathbb{C}^{k+1}$.

   Note that $L\subset \mathbb{C}P^{k}\times \mathbb{C}^{k+1}$ is
   monotone since $\mathbb{C}P^{k}$ is a monotone symplectic manifold
   and $H_{1}(L,\mathbb{Z})=0$. A simple computation shows that the
   minimal Maslov number of $L$ is $N_{L}=2k+2$. Thus $N_{L} \geqslant
   2$ and the Floer homology of $L$ is well defined.

   We now claim that $HF(L,L)=0$. Indeed, L can be disjoint from
   itself by a large enough linear translation along the
   $\mathbb{C}^{k+1}$-factor, and linear translations are
   Hamiltonian.

   We now compare this to a computation using Oh's spectral sequence.
   We use here the notations from Section 2. Since $N_{L}=2k+2$ and
   $\dim L=2k+1$ we must have $\partial_{j}=0 \quad \forall \,
   j\geqslant 2$ , hence $d_{F}=\partial_{0}+\partial_{1}$.  Therefore
   $H^*([\partial_{1}])=HF(L,L)$. Since $HF(L,L)=0$ this implies that
   $H^{i+2k+1}(L,\mathbb{Z}_{2}) \xrightarrow{[\partial_{1}]} H^i(L;
   \mathbb{Z}_2) \xrightarrow{[\partial_{1}]} H^{i-2k-1}(L;
   \mathbb{Z}_2)$ is an exact sequence $\forall \, 0\leqslant
   i\leqslant 2k+1$.  As $\dim L=2k+1$ we immediately get that
   $H^{i}(L,\mathbb{Z}_{2})=0 \quad \forall \, 1\leqslant i\leqslant
   2k$ and
   $H^{0}(L,\mathbb{Z}_{2})=H^{2k+1}(L,\mathbb{Z}_{2})=\mathbb{Z}_{2}$,
   or in other words that $H^{*}(L,\mathbb{Z}_{2})\cong
   H^{*}(S^{2k+1},\mathbb{Z}_{2})$.
\end{proof}

%
%
\begin{proof}[Proof of Theorem~\ref{T:lens}]
   Denote by $i : L \hookrightarrow T^*(Lens_{m}^{2k+1})$ the
   inclusion. We claim that $i_{*}:\pi_{1}(L)\rightarrow
   \pi_{1}(T^{*}(Lens_{m}^{2k+1})) \cong \mathbb{Z}_m$ is an
   isomorphism.  Indeed, suppose on the contrary that $i_*: \pi_1(L)
   \to \pi_1(T^*(Lens_{m}^{2k+1}))$ is not surjective. Let $Y \to
   Lens_{m}^{2k+1}$ be the covering associated to the subgroup
   $i_*(\pi_1(L)) \subset \pi_1(T^*(Lens_{m}^{2k+1})) \cong
   \pi_1(Lens_{m}^{2k+1})$. Consider the corresponding covering
   $T^*(Y) \to T^*(Lens_{m}^{2k+1})$, and let $j: L \hookrightarrow
   T^*(Y)$ be a lifting of $i:L \hookrightarrow T^*(Lens_{m}^{2k+1})$
   (it exists because Y is chosen to be a covering associated to the
   subgroup $i_*(\pi_1(L)) \subset \pi_1(T^*(Lens_{m}^{2k+1})) \cong
   \pi_1(Lens_{m}^{2k+1})$, and it is obviously a Lagrangian embedding).
   Since $i_*$ is not surjective there exists a nontrivial deck
   transformation $\sigma: Y \to Y$. We claim that $\sigma$ must
   be isotopic to the identity. To see this, note that all connected
   coverings of $Lens_{m}^{2k+1}$ are determined, up to an isomorphism,
   by a subgroup  $ \mathbb{Z}_{d} \hookrightarrow \pi_{1}(Lens_{m}^{2k+1})
   \cong \mathbb{Z}_{m}$. Hence every such covering is isomorphic
   to $ Lens_{d}^{2k+1} \rightarrow Lens_{m}^{2k+1} $, where
   $ d | m $. But for such coverings it is easy to check that all
   deck transformations are isotopic to the identity.
   Thus $\sigma$ is isotopic to the identity, and therefore
   the canonical lift of $\sigma$, $\Phi_{\sigma}: T^*(Y)
   \to T^*(Y)$ is Hamiltonian. Now $\Phi_{\sigma}$ itself is a
   nontrivial deck transformation for the covering $T^*(Y) \to
   T^*(Lens_{m}^{2k+1})$. Therefore since $j: L \hookrightarrow
   T^*(Y)$ is a lifting of $i: L \hookrightarrow T^*(Lens_{m}^{2k+1})$
   we must have $\Phi_{\sigma}(j(L)) \cap j(L) = \emptyset$. But this
   is impossible since $j(L) \subset T^*(Y)$ is an exact Lagrangian.
   We therefore conclude that $i_*: \pi_1(L) \to
   \pi_1(T^*(Lens_{m}^{2k+1}))$ is surjective and therefore an
   isomorphism.

   Consider now the universal coverings $\widetilde{L} \to L$ and
   $T^*(S^{2k+1}) \to T^*(Lens_{m}^{2k+1})$. As $i_*$ is an
   isomorphism we can lift $\widetilde{L}$ into $T^*(S^{2k+1})$ namely
   we have a commutative diagram
   \begin{equation*}
      \begin{array}{clcr}
         \widetilde{L} & \hookrightarrow & T^{*}(S^{2k+1}) \\
         \downarrow &   &  \downarrow \\
         L & \hookrightarrow & T^{*}(Lens_{m}^{2k+1})
      \end{array}
   \end{equation*}
   where the upper horizontal map is a Lagrangian embedding. By
   Theorem~\ref{T:sphere} we have $H^{*}(\widetilde{L}; \mathbb{Z}_2)
   \cong H^{*}(S^{2k+1}; \mathbb{Z}_2)$.

   Now we use the Cartan-Leray spectral sequence for coverings (see
   e.g.~\cite{Book-Spectral}), in order to derive $ H^*(L,
   \mathbb{Z}_2) $ from $ H^{*}(\widetilde{L}, \mathbb{Z}_2) $.
   Let us outline the main steps in this computation.

   Recall that given a covering  $X \rightarrow Y$ with group $\pi$
   and a commutative group $ A $, there exists a spectral
   sequence with $E_{p,q}^{2}=H_{p}(\pi,H_{q}(X,A))$ that
   converges to $H_{*}(Y,A)$. Here $H_{p}(\pi,H_{q}(X,A))$
   stands for the homology of the group $\pi$  with coefficients
   in the $\pi$-module $H_{q}(X,A)$.
   (Namely, $H_{p}(\pi,H_{q}(X,A)) =
   Tor_{p}^{\mathbb{Z}[\pi]}(\mathbb{Z},H_{q}(X,A))$.)

   In our situation we have $ A = \mathbb{Z}_{2} $, and the
   covering $ \widetilde{L} \rightarrow L$ with group $\pi =
   \mathbb{Z}_{m} $. Since $ H_{*}(\widetilde{L},\mathbb{Z}_{2})
   \cong H_{*}(S^{2k+1},\mathbb{Z}_{2}) $ and since $\pi$ acts
   trivially on $ H_{0}(\widetilde{L},\mathbb{Z}_{2}) $ and on
   $ H_{2k+1}(\widetilde{L},\mathbb{Z}_{2}) $, the beginning of
   the spectral sequence looks as follows:


\[  \begin{array}{c|c|c|c|c|c|c|c|r}
         \vdots  &  \vdots &  \vdots &  \vdots &  \vdots & \vdots & \vdots & \vdots & \\  \cline{2-9}
       2k+2 & 0 & 0 & \ldots & 0 & 0 & 0 & \cdots \\ \cline{2-9}
       2k+1 & H_{0}(\mathbb{Z}_{m},\mathbb{Z}_{2}) & H_{1}(\mathbb{Z}_{m},\mathbb{Z}_{2}) & \cdots &
       H_{2k}(\mathbb{Z}_{m},\mathbb{Z}_{2}) & H_{2k+1}(\mathbb{Z}_{m},\mathbb{Z}_{2}) &
       0 & \cdots
       \\ \cline{2-9}
       2k & 0 & 0 & \cdots & 0 & 0 & 0 & \cdots \\  \cline{2-9}
       \vdots & \vdots & \vdots & \ddots & \vdots & \vdots & \vdots & \cdots \\  \cline{2-9}
       1 & 0 & 0 & \cdots & 0 & 0 & 0 & \cdots \\  \cline{2-9}
       0 & H_{0}(\mathbb{Z}_{m},\mathbb{Z}_{2}) & H_{1}(\mathbb{Z}_{m},\mathbb{Z}_{2}) & \cdots &
       H_{2k}(\mathbb{Z}_{m},\mathbb{Z}_{2}) & H_{2k+1}(\mathbb{Z}_{m},\mathbb{Z}_{2}) &
       0 & \cdots
       \\  \cline{2-9}

    \end{array}
\]

$ \hspace{22 mm}  0   \hspace{22 mm}  1 \hspace{11.5 mm} \cdots
\hspace{12 mm}  2k \hspace{20 mm}  2k+1 \hspace{2 mm}  $

   Since the differential at the $ r $'th step is $ d_{r} :
   E_{p,q}^{r} \rightarrow E_{p-r,q+r-1}^{r} $ we get
   $ H_{i}(L,\mathbb{Z}_{2}) =
   \bigoplus_{p+q=i}E_{p,q}^{\infty} \cong
   H_{i}(\mathbb{Z}_{m},\mathbb{Z}_{2})$
   for every $ 0 \leqslant i \leqslant 2k $, and of course $ H_{2k+1}(L,\mathbb{Z}_{2})
   = \mathbb{Z}_{2} $.
   A straightforward computation of $
   H_{i}(\mathbb{Z}_{m},\mathbb{Z}_{2}) $ now gives:

   \[H_{i}(L,\mathbb{Z}_{2})=
   \left\{\begin{array}{ll}
            \mathbb{Z}_{2} & \mbox{$m$ even} \\
            0 & \mbox{$m$ odd} \\
          \end{array}
   \right. \]
   for every $ 0 \leqslant i \leqslant 2k $ .
   On the other hand, this is precisely the homology $
   H_{i}(Lens_{m}^{2k+1},\mathbb{Z}_{2}) $.

\end{proof}

\subsubsection*{Acknowledgments}
I would like to thank my supervisor Paul Biran for the help and
attention he gave me.

\end{document}